\documentclass[preprint,3p,number,sort&compress]{elsarticle}

\usepackage{amsmath,amssymb,calrsfs}
\usepackage[hidelinks]{hyperref}

\newtheorem{theorem}{Theorem}[section]
\newtheorem{lemma}[theorem]{Lemma}
\newtheorem{definition}[theorem]{Definition}
\newtheorem{corollary}[theorem]{Corollary}
\newtheorem{remark}[theorem]{Remark}
\newproof{proof}{Proof}

\newcommand{\proofofref}{}
\newproof{zproofof}{Proof of \proofofref}
\newenvironment{proofof}[1]
 {\renewcommand{\proofofref}{#1}\zproofof}
 {\endzproofof}

\numberwithin{equation}{section}
\newcommand{\D}{{\scriptscriptstyle\Delta}}

\begin{document}
\begin{frontmatter}
\title{Self-improving boundedness\\ of the maximal operator on quasi-Banach lattices\\ over spaces of homogeneous type}

\author{Alina Shalukhina}  
\ead{a.shalukhina@campus.fct.unl.pt}

\address{
Centro de Matem\'atica e Aplica\c{c}\~oes,
Departamento de Matem\'a\-tica,
Faculdade de Ci\^en\-cias e Tecnologia,\\
Universidade NOVA de Lisboa,
Quinta da Torre,
2829-516 Caparica, 
Portugal}

\begin{abstract}
We prove the self-improvement property of the Hardy--Littlewood maximal operator on quasi-Banach lattices with the 
Fatou property in the setting of spaces of homogeneous type. Our result is a generalization of the boundedness criterion
obtained in 2010 by Lerner and Ombrosi for maximal operators on quasi-Banach function spaces over Euclidean spaces. The
specialty of the proof for spaces of homogeneous type lies in using adjacent grids of Hyt\"onen--Kairema dyadic 
cubes and studying the maximal operator alongside its ``dyadic'' version. Then we apply the obtained result to 
variable Lebesgue spaces over spaces of homogeneous type.
\end{abstract}

\end{frontmatter}

\section{Introduction}
The primary interest in the study of the Hardy--Littlewood maximal operator, classically defined for a measurable function 
$f$ on $\mathbb{R}^n$ by 
\[
Mf(x)=\sup_{B\ni x}\frac1{|B|}\int_{B} |f(y)|\,dy,
\]
where the supremum of the integral means of $f$ is taken over all balls $B$ containing a point $x\in\mathbb{R}^n$, has
always been related to finding necessary and sufficient conditions for the operator's boundedness on various function spaces. 
Over the last two decades, characterization of spaces on which $M$ is bounded became a subject of extensive research in
the setting of variable Lebesgue spaces $L^{p(\cdot)}(\mathbb{R}^n)$. This endeavor has recently culminated in
Lerner's criterion for the boundedness of $M$ on $L^{p(\cdot)}(\mathbb{R}^n)$ forthcoming in~\cite{L24}.
Along with the effort to obtain such a characterization in terms of the exponent function, another natural question was
raised: if $M$ is bounded on a certain $L^{p(\cdot)}$ space, which perturbations of the exponent $p(\cdot)$ preserve 
the boundedness? Once we have the boundedness of $M$ on one space, can it be automatically extended to a family 
of the closely related spaces?

A significant development in this direction was the 2005 result of Diening \cite[Theorem~8.1]{D05} within which he 
proved, in particular, that $M$ is bounded on a space $L^{p(\cdot)}(\mathbb{R}^n)$ with a bounded exponent 
$p(\cdot)$ if and only if it is bounded on $L^{p(\cdot)/q}(\mathbb{R}^d)$ for some $q>1$ (see also 
\cite[Theorem~5.7.2]{DHHR11}). This property received the name of ``left-openness'' in his work---the terminology
borrowed from the theory of Muckenhoupt weights, whose classes are ``left-open'' in the sence that for a weight 
$w\in A_p$, there exists $q<p$ such that $w\in A_q$. Alternatively, this classical feature of the Muckenhoupt classes 
may be referred to as a self-improving property, so the boundedness result for the maximal operator adopts this name 
later as well. 

Two years after Diening's result, Lerner and P\'erez~\cite[Corollary~1.3]{LP07} proved the self-improving boundedness 
of $M$ in the more general setting of quasi-Banach function spaces $X(\mathbb{R}^n)$. Thus this property, known for
weighted Lebesgue spaces (as a consequence of the left-openness of the Muckenhoupt classes) and variable 
$L^{p(\cdot)}$ spaces (due to Diening), received a generalization and a unified proof. Expressed in the equivalent terms 
of the parameterized maximal function $M_rf=M(|f|^r)^{1/r}$ instead of ``perturbations'' of the space, the generalized
result stated that $M$ is bounded on $X$ if and only if $M_r$ is bounded on $X$ for some $r>1$. 

Subsequently, in 2010, Lerner and Ombrosi \cite[Theorem~1.2]{LO10} complemented the self-improving property on 
quasi-Banach spaces $X$ over $\mathbb{R}^n$  by considering more general maximal operators $M_{\mathcal{B}}$ 
with respect to a basis $\mathcal{B}$ and giving an additional ``self-improvement'' characterization in terms of
$M_{\mathcal{B},r}$ with $r<1$. It is exactly this result that we took as a reference point for our work---and transferred
it, for the classical maximal operator, to a different topological setting. Except for switching from the abstract
$M_{\mathcal{B}}$ to the classical $M$, our main Theorem~\ref{main-theorem} can be viewed as a generalization of the
Lerner--Ombrosi result: on the topological level, we replace the Euclidean space $\mathbb{R}^n$ by a space of
homogeneous type $(\Omega,d,\mu)$, which is essentially a quasi-metric measure space with a doubling measure; 
on the functional level, we notice that the requirement on $X$ to be a quasi-Banach function space is abundant for the
proof---the framework of which we follow, though bringing in techniques specific for the new topology---and so assume 
only that $X$ is a quasi-Banach lattice with the Fatou property (see~\cite{LN23} for differences between the above 
two concepts).    

We begin in Section~\ref{sec:lattices} by defining a quasi-Banach lattice and showing that any variable Lebesgue
space with the exponent function bounded away from zero falls within its definition: this tangible example gives better
perception of the abstract reasoning and is a natural setting for the later application of the main result. To make the
connection between our generalization of the Lerner--Ombrosi theorem and the original self-improvement result for
$L^{p(\cdot)}$ even more intuitive, we formulate Theorem~\ref{main-theorem} in terms of convexifications $X^{(r)}$ 
of a lattice $X$. These are logically the same ``perturbations'' of the space as $L^{rp(\cdot)}$ would be with respect 
to $L^{p(\cdot)}$, and we provide their precise definition in Subsection~\ref{subsec:convex} following
Maligranda~\cite[p.~99]{M04}.

Having said that, we now present the main theorem.
\begin{theorem}\label{main-theorem}
Let $X(\Omega,d,\mu)$ be a quasi-Banach lattice over a space of homogeneous type $(\Omega,d,\mu)$. If 
$X(\Omega,d,\mu)$ has the Fatou property, the following statements are equivalent: 
\begin{enumerate}
\item[\textnormal{(1)}] $M$ is bounded on $X(\Omega,d,\mu)$.

\item[\textnormal{(2)}] For all $s>1$, $M$ is bounded on $X^{(s)}(\Omega,d,\mu)$ and 
\[
\lim_{s\to1^+}(s-1)\|M\|_{X^{(s)}\to X^{(s)}}=0.
\]

\item[\textnormal{(3)}] There exists $r_0\in(0,1)$ such that if $r\in[r_0,1)$, then $M$ is bounded on 
$X^{(r)}(\Omega,d,\mu)$.
\end{enumerate}
\end{theorem}
Technically, the above collection of equivalent statements splits into the two ``trivial'' implications $(1)\Rightarrow(2)$ 
and $(3)\Rightarrow(1)$ and the converse non-trivial implications $(2)\Rightarrow(1)$ and $(1)\Rightarrow(3)$ containing
the self-improvement property itself. The former are a simple consequence of H\"older's inequality and extend the
boundedness of $M$ to ``higher'' convexifications of the lattice $X$; the latter provide the extension to a range of 
``lower'' convexifications and require a sort of reverse-H\"older bound for the maximal operator. 

In our setting of spaces of homogeneous type---the review of which comprises the entire Section~\ref{sec:SHT}---such 
a bound is available for the so-called ``dyadic'' maximal function $M^\mathcal{D}w$ of a ``dyadic'' $A_1$ weight $w$. 
By the ``dyadic'' objects, introduced in Subsections~\ref{subsec:MOs} and \ref{subsec:A1-weights}, we always 
understand those defined through the adjacent grids of Hyt\"onen--Kairema dyadic cubes (see \cite[Section~4]{HK12})
instead of balls. Then the Rubio de Francia iteration algorithm, which we apply after Lerner and Ombrosi, allows 
construction of appropriate ``dyadic'' $A_1$ weights for the use of  the reverse-H\"older type bound
(Corollary~\ref{cor:point-est}) based on the weak reverse H\"older inequality originally proved for the ``dyadic'' $A_\infty$
weights by Anderson, Hyt\"onen and Tapiola~\cite[Theorem~5.4]{AHT17}. Given also the equivalence between
$M^\mathcal{D}$ and $M$, we are able to switch between the two operators when convenient and eventually prove 
the non-trivial implications using this trick.

Section~\ref{sec:proof} presents the proof of Theorem~\ref{main-theorem} in line with its division into separate
implications---with Lemma~\ref{le:trivial-impl} uniting the trivial implications and Theorems~\ref{th:main-lower} and
\ref{th:main-upper} corresponding to the ``self-improvement'' part. A great help in organizing our proof was 
Theorem~4.37 from the book by Cruz-Uribe and Fiorenza~\cite{CF13}, which contains an adaptation of the Lerner--Ombrosi
theorem to the special case of spaces $L^{p(\cdot)}(\mathbb{R}^n)$ with the ``Banach'' range of exponents $p(\cdot)$
not going below 1. Towards the close, in Subsection~\ref{subsec:appl}, we give an application of the main result to the
variable Lebesgue spaces $L^{p(\cdot)}(\Omega,d,\mu)$ with $p(\cdot)$ bounded away from zero. 

Let us also note that the reference paper by Lerner and Ombrosi~\cite{LO10} set in quasi-Banach spaces has recently 
given rise to a few further studies of the self-improvement property for maximal operators on Banach function spaces. 
Thus, e.g., a sharp version of the Lerner--Ombrosi theorem appeared in the work of Nieraeth~\cite[Theorem~2.34]{N23} 
for the maximal operator with respect to some general bases in the setting of Banach function spaces. Another similar
self-improvement property of $M$ on $r$-convex Banach function spaces was proved by Lorist and
Nieraeth~\cite[Theorem~3.1]{LN24}, with the possibility to extend this result to Banach function spaces over spaces
of homogeneous type~\cite[Remark~3.5]{LN24}.

\section{Quasi-Banach Lattices and Their Convexifications}\label{sec:lattices}

\subsection{Quasi-Banach lattices}
We denote by $L^0(\Omega,\mu)$ the space of measurable complex-valued functions on a measure space 
$(\Omega,\mu)$. Let $X(\Omega,\mu)\subset L^0(\Omega,\mu)$ be a quasi-normed space---that is, a linear subspace of 
$L^0(\Omega,\mu)$ endowed with a quasi-norm $\|\cdot\|_X$, which differs from the usual norm by the weakened
triangle inequality
\[
\|f+g\|_X\le C_\D (\|f\|_X+\|g\|_X)
\]
holding with a constant $C_\D\ge1$ independent of $f$ and $g$. The space $X(\Omega,\mu)$ is called a quasi-normed
lattice if it additionally satisfies the following \textit{lattice property}:
\begin{center}
if $f\in X(\Omega,\mu)$, and $g\in L^0(\Omega,\mu)$ is such that $|g|\le|f|$,\\
then $g\in X(\Omega,\mu)$ and $\|g\|_X\le\|f\|_X$.
\end{center}

Quasi-normed lattices and their complete counterparts, quasi-Banach lattices, form the functional setting for our work.
In the main result of this paper, a quasi-Banach lattice $X(\Omega,\mu)$, though taken over a specific measure space
which is a space of homogeneous type, will also be required to have a property stronger than completeness---the 
\textit{Fatou property} with constant $C_\mathcal{F}>0$: 
\begin{center}
if $0\le f_n\uparrow f$ for a sequence $\{f_n\}\subset X(\Omega,\mu)$ and $\sup_{n\ge0}\|f_n\|_X<\infty$,\\
then $f\in X(\Omega,\mu)$ and $\|f\|_X\le C_\mathcal{F}\sup_{n\ge0}\|f_n\|_X$.
\end{center}
The proof that the Fatou property of a quasi-normed lattice indeed imples its completeness is outlined
in~\cite[Remark~2.1(ii)]{LN23}; although a stronger form of the Fatou property is assumed there, the same argument
applies unchanged within our definition. (For the case of normed lattices, see also \cite[Ch.~IV, \S~3, Theorem~4]{KA82}.)

For quasi-Banach lattices with the Fatou property, we have a useful version of the Aoki-Rolewicz theorem 
for infinite sums.
\begin{theorem}\label{th:Aoki-Rol-inf}
Let a quasi-Banach lattice $X(\Omega,\mu)$ have the Fatou property with constant $C_\mathcal{F}$. Then 
for any nonnegative sequence $\{f_k\}\subset X(\Omega,\mu)$, there holds
\[
\left\|\sum_{k=0}^\infty f_k\right\|_X
\le 
2^{1/\rho} C_\mathcal{F}\left(\sum_{k=0}^\infty \|f_k\|_X^\rho\right)^{1/\rho},
\]
where the number $\rho\in(0,1]$ is given by $2^{{1/\rho}-1}=C_\D$ and will hereafter be referred to as 
the Aoki-Rolewicz exponent of the lattice $X(\Omega,\mu)$.
\end{theorem}
\begin{proof}
For $n=0,1,\ldots$, let us denote $F_n=\sum_{k=0}^n f_k$ and $F=\sum_{k=1}^\infty f_k$; then $F_n\uparrow F$. 
If $\sum_{k=0}^\infty \|f_k\|_X^\rho=\infty$, the statement is trivial. If
\[
\sum_{k=0}^\infty \|f_k\|_X^\rho=:M<\infty, 
\]
then for all $n$, we have $\sum_{k=0}^n \|f_k\|_X^\rho \le M$.
By the Aoki-Rolewicz theorem for finite sums (see~\cite[p.~47]{K86}; cf. a weaker 
formulation in~\cite[p.~3]{KPR84} with $4^{1/\rho}$ in place of $2^{1/\rho}$),
for every $n\ge0$,
\[
\|F_n\|_X \le 2^{1/\rho}\left(\sum_{k=0}^n \|f_k\|_X^\rho\right)^{1/\rho} \le 2^{1/\rho} M^{1/\rho}<\infty,
\]
hence $\sup_{n\ge0}\|F_n\|_X<\infty$. The Fatou property implies that $F\in X(\Omega,\mu)$ and
\[
\|F\|_X \le C_\mathcal{F} \sup_{n\ge0}\|F_n\|_X
\le 2^{1/\rho} C_\mathcal{F} M^{1/\rho} 
= 2^{1/\rho} C_\mathcal{F} \left(\sum_{k=0}^\infty \|f_k\|_X^\rho\right)^{1/\rho},
\]
which is the desired inequality.
\end{proof}

\subsection{Convexifications}\label{subsec:convex}
Together with quasi-normed lattices $X(\Omega,\mu)$, we consider their $r$-convexifications, $r>0$, defined by
\[
X^{(r)}(\Omega,\mu):=\{f\in L^0(\Omega,\mu)\ :\ |f|^r\in X(\Omega,\mu)\}.
\]
It is well-known that each of the convexifications is a quasi-normed lattice itself. This fact is mentioned without proof, e.g., 
in~\cite[p.~99]{M04}; being unable to provide a better reference, we prove it here.
\begin{lemma}\label{le:convex-lattice}
If $X(\Omega,\mu)$ is a quasi-normed lattice with a quasi-norm $\|\cdot\|_X$, then $X^{(r)}(\Omega,\mu)$, for
any $r>0$, is a quasi-normed lattice with the quasi-norm 
\[
\|f\|_{X^{(r)}}:=\|\,|f|^r \|_X^{1/r}.
\]
\end{lemma}
\begin{proof}
Fix $r>0$, and let us check first that $X^{(r)}(\Omega,\mu)$ is a vector space. Take two functions 
$f,g\in X^{(r)}(\Omega,\mu)$. Since there holds the pointwise relation
\begin{equation}\label{eq:estimates}
|f+g|^r\le(|f|+|g|)^r\le\max\{1,2^{r-1}\}(|f|^r+|g|^r),
\end{equation}
and $|f|^r,|g|^r\in X(\Omega,\mu)$ by our choice, the function on the right-hand side of the above inequality belongs
to $X(\Omega,\mu)$, which is a vector space. Then, by the lattice property, the smaller left-hand function 
$|f+g|^r\in X(\Omega,\mu)$, and thus $f+g\in X^{(r)}(\Omega,\mu)$. Trivially, we have 
$|\lambda f|^r=|\lambda|^r|f|^r\in X(\Omega,\mu)$, which implies that $\lambda f\in X^{(r)}(\Omega,\mu)$, for any 
$\lambda\in\mathbb{C}$. Thus, an $r$-convexification of a lattice is indeed a vector space.

Next, it follows from~\eqref{eq:estimates} and the quasi-triangle inequality for $X(\Omega,\mu)$, which holds with the
``quasi-norm'' constant $C_\D(X)$, that
\begin{align*}
\|f+g\|_{X^{(r)}}&=\|\,|f+g|^r\|_X^{1/r} \le (\max\{1,2^{r-1}\})^{1/r} \|\,|f|^r+|g|^r\|_X^{1/r} \\
&\le\max\{1,2^{1-1/r}\}\, C_\D(X)^{1/r} (\|\,|f|^r\|_X+\|\,|g|^r\|_X)^{1/r} \\
&\le\max\{1,2^{1-1/r}\}\, C_\D(X)^{1/r} \max\{1,2^{1/r-1}\} (\|\,|f|^r\|_X^{1/r}+\|\,|g|^r\|_X^{1/r}) \\
&=C_\D(X^{(r)}) (\|f\|_{X^{(r)}}+\|g\|_{X^{(r)}}),
\end{align*}
where 
$C_\D(X^{(r)}):=2^{|1-1/r|}C_\D(X)^{1/r}$. Hence, $\|\cdot\|_{X^{(r)}}$ satisfies the quasi-triangle inequality with the constant $C_\D(X^{(r)})$, about which we note, for later use, that 
\begin{equation}\label{eq:C-tri-convex}
C_\D(X^{(r)})\le2C_\D(X)\,\text{ whenever }\,r\ge1.
\end{equation} 
Verifying the other axioms of the quasi-norm for $\|\cdot\|_{X^{(r)}}$ and establishing the lattice property in the space 
$X^{(r)}(\Omega,\mu)$ are straightforward.
\end{proof}

Moreover, the action of ``convexifying'' a lattice preserves the Fatou property, and hence completeness.
\begin{lemma}\label{le:Fatou}
If a quasi-Banach lattice $X(\Omega,\mu)$ has the Fatou property with constant $C_\mathcal{F}(X)$, then its
\mbox{$r$-convexification} $X^{(r)}(\Omega,\mu)$, $r>0$, is a quasi-Banach lattice satisfying the Fatou property with 
the constant
\[
C_\mathcal{F}(X^{(r)})=C_\mathcal{F}(X)^{1/r}.
\]
\end{lemma}
\begin{proof}
Let us take a sequence $\{f_n\}$ in $X^{(r)}(\Omega,\mu)$ such that $0\le f_n\uparrow f$ and 
\[
\sup_{n\ge0}\|f_n\|_{X^{(r)}}=\sup_{n\ge0}\|f_n^r\|_X^{1/r}<\infty.
\]
Applying the Fatou property to the sequence $0\le f_n^r\uparrow f^r$ in $X(\Omega,\mu)$, we conclude that 
$f^r\in X(\Omega,\mu)$---or equivalently, $f\in X^{(r)}(\Omega,\mu)$---and
\[
\|f\|_{X^{(r)}}=\|f^r\|_X^{1/r}\le\Big[C_\mathcal{F}(X) \sup_{n\ge0}\|f_n^r\|_X\Big]^{1/r}
=C_\mathcal{F}(X)^{1/r} \sup_{n\ge0}\|f_n\|_{X^{(r)}}.
\]
Thus, $X^{(r)}(\Omega,\mu)$ satisfies the Fatou property with the constant $C_\mathcal{F}(X)^{1/r}$.
\end{proof}

\subsection{Example: variable Lebesgue spaces}\label{subsec:example}
For a function $p(\cdot):\Omega\to(0,\infty]$ measurable on $(\Omega,\mu)$, called an exponent function, and an
arbitrary $f\in L^0 (\Omega,\mu)$, consider the modular functional associated with $p(\cdot)$ given by
\[
m_{p(\cdot)}(f):=\int_{\Omega\setminus\Omega_\infty}|f(x)|^{p(x)}d\mu(x)+
\operatornamewithlimits{ess\,sup}\limits_{x\in\Omega_\infty}|f(x)|, 
\]
where $\Omega_\infty=\{x\in\Omega\,:\,p(x)=\infty\}$. By customary definition, the variable Lebesgue space
$L^{p(\cdot)}(\Omega,\mu)$ consists of all measurable functions $f$ such that $m_{p(\cdot)}(f/\lambda)<\infty$ 
for some $\lambda>0$ depending on $f$. 

When an exponent function $p(\cdot)$ is essentially bounded away from zero, the space $L^{p(\cdot)}(\Omega,\mu)$ 
is a quasi-Banach space---for the setting of $\mathbb{R}^n$, this fact is mentioned without proof in~\cite[p.~940]{KV14}.
More generally, the following statement is true. 
\begin{theorem}\label{th:Q-B-lattice}
Given an exponent function $p(\cdot)$ such that 
\[
p_-:=\operatornamewithlimits{ess\,inf}\limits_{x\in\Omega}p(x)>0,
\]
the space $L^{p(\cdot)}(\Omega,\mu)$ is a quasi-Banach lattice, possessing the Fatou property, with respect to 
the Luxemburg--Nakano quasi-norm 
\begin{equation}\label{eq:Luxemb-QN}
\|f\|_{p(\cdot)}:=\inf\{\lambda>0\ :\ m_{p(\cdot)}(f/\lambda)\le1\}.
\end{equation}
Further, the ``quasi-norm'' constant is $C_\D=\max\{1,2^{1/p_- -1}\}$.
\end{theorem}

The proof of Theorem~\ref{th:Q-B-lattice} employs the next elementary properties of the modular.
\begin{lemma}\label{le:elem-prop}
For any exponent function $p(\cdot)$, the following are true:
\begin{enumerate}
\item[\textnormal{(i)}] $m_{p(\cdot)}$ is order preserving: if $|g|\le|f|$, then 
$m_{p(\cdot)}(g)\le m_{p(\cdot)}(f)$.

\item[\textnormal{(ii)}] For any $0<\alpha<1$, 
\begin{equation}\label{const-in-modular}
m_{p(\cdot)}(\alpha f)\le\alpha^{\min\{p_-,1\}} m_{p(\cdot)}(f).
\end{equation}

\item[\textnormal{(iii)}] If $p_-<1$, then given $\alpha,\beta\ge0$ such that $\alpha+\beta=1$, there holds 
\begin{equation}\label{eq:kind-of-convexity}
m_{p(\cdot)}(\alpha f+\beta g)\le\alpha^{p_-}m_{p(\cdot)}(f)+\beta^{p_-}m_{p(\cdot)}(g).
\end{equation}
\end{enumerate}
\end{lemma}
\begin{proof}
Properties~(i) and (ii) follow simply from the definition of the modular. For~(iii), observe first that if $p_-<1$, then 
for any $x\in\Omega$ and each $y,z\ge0$, we have 
\begin{equation}\label{eq:aux-ineq}
(\alpha y+\beta z)^{p(x)}\le\alpha^{p_-}y^{p(x)}+\beta^{p_-}z^{p(x)}.
\end{equation} 
Indeed, if $p(x)\ge1$, then by Jensen's inequality, 
\[
(\alpha y+\beta z)^{p(x)}\le \alpha y^{p(x)}+\beta z^{p(x)}\le \alpha^{p_-}y^{p(x)}+\beta^{p_-}z^{p(x)};
\]
otherwise, if $0<p(x)<1$, we use a simple inequality from~\cite[p.~121]{CR16} to deduce
\[
(\alpha y+\beta z)^{p(x)}\le (\alpha y)^{p(x)}+(\beta z)^{p(x)}\le \alpha^{p_-}y^{p(x)}+\beta^{p_-}z^{p(x)}.
\]

Applying~\eqref{eq:aux-ineq}, we conclude that
\begin{align*}
m_{p(\cdot)}(\alpha f+\beta g)&\le 
\alpha^{p_-}\int_{\Omega\setminus\Omega_\infty}|f(x)|^{p(x)}d\mu(x)
+ \beta^{p_-}\int_{\Omega\setminus\Omega_\infty}|g(x)|^{p(x)}d\mu(x) \\
&\quad+\alpha\operatornamewithlimits{ess\,sup}_{x\in\Omega_\infty}|f(x)|
+\beta\operatornamewithlimits{ess\,sup}_{x\in\Omega_\infty}|g(x)| \\
&\le\alpha^{p_-}m_{p(\cdot)}(f)+\beta^{p_-}m_{p(\cdot)}(g),
\end{align*}
which is exactly the desired inequality~\eqref{eq:kind-of-convexity}.
\end{proof}
\begin{proofof}{Theorem~\ref{th:Q-B-lattice}}
We begin by checking that \eqref{eq:Luxemb-QN} indeed defines a quasi-norm on $L^{p(\cdot)}(\Omega,\mu)$, 
namely, that the functional $\|\cdot\|_{p(\cdot)}$ has the following properties:
\begin{enumerate}
\item[(a)] $\|f\|_{p(\cdot)}=0$ if and only if $f=0$;

\item[(b)] $\|\alpha f\|_{p(\cdot)}=|\alpha|\|f\|_{p(\cdot)}$ for all $\alpha\in\mathbb{C}$;

\item[(c)] $\|f+g\|_{p(\cdot)}\le C_\D(\|f\|_{p(\cdot)}+\|g\|_{p(\cdot)})$, where $C_\D=\max\{1,2^{1/p_- -1}\}$.
\end{enumerate}

Clearly, $\|0\|_{p(\cdot)}=0$. Suppose $\|f\|_{p(\cdot)}=0$; then $m_{p(\cdot)}(f/\lambda)\le1$ for all $\lambda>0$.
We immediately have that $|f(x)|\le\lambda$ for almost every $x\in\Omega_\infty$, hence $f=0$ 
on $\Omega_\infty$. At the same time, if $0<\lambda<1$, then 
\[
1\ge\int_{\Omega\setminus\Omega_\infty} \left|\frac{f(x)}{\lambda}\right|^{p(x)}d\mu(x)\ge
\lambda^{-p_-}\int_{\Omega\setminus\Omega_\infty}|f(x)|^{p(x)}d\mu(x),
\]
which implies
\[
\int_{\Omega\setminus\Omega_\infty}|f(x)|^{p(x)}d\mu(x)\le\lim_{\lambda\to0^+}\lambda^{p_-}=0,
\]
and therefore $f=0$ on $\Omega\setminus\Omega_\infty$. Thus $f=0$ and we have proved (a). 

For (b), note that if $\alpha=0$, this follows from (a). Otherwise, if $\alpha\ne0$, we get
\begin{align*}
\|\alpha f\|_{p(\cdot)}=\inf&\left\{\lambda>0\ :\ m_{p(\cdot)}\left(\frac{|\alpha|f}{\lambda}\right)\le1\right\} \\
&\quad=|\alpha|\inf\left\{\frac{\lambda}{|\alpha|}>0\ :\ m_{p(\cdot)}\left(\frac{f}{\lambda/|\alpha|}\right)\le1\right\}
=|\alpha|\|f\|_{p(\cdot)}.
\end{align*}

As to (c), it is well known that $\|\cdot\|_{p(\cdot)}$ defines a norm on $L^{p(\cdot)}$ when $p_-\ge1$, so the
usual triangle inequality holds in this case (see, e.g., \cite[Theorem~2.17]{CF13} for the proof in $\mathbb{R}^n$). 
Therefore, it remains to establish the quasi-triangle inequality with the constant $C_\D=2^{1/p_- -1}$ when $0<p_-<1$. 

Fix $\lambda_f>\|f\|_{p(\cdot)}$ and $\lambda_g>\|g\|_{p(\cdot)}$; then $m_{p(\cdot)}(f/\lambda_f)\le1$ and 
$m_{p(\cdot)}(g/\lambda_g)\le1$ due to the order preserving property. Let $\lambda=\lambda_f+\lambda_g$. Using
inequality~\eqref{eq:kind-of-convexity} and the fact that $t\mapsto t^{p_-}$ is concave and hence
\[
y^{p_-}+z^{p_-}\le2^{1-p_-}(y+z)^{p_-}\ \text{ for each }y,z\ge0,
\]
we find that
\begin{align*}
m_{p(\cdot)}\left(\frac{f+g}{\lambda}\right)&=
m_{p(\cdot)}\left(\frac{\lambda_f}{\lambda}\frac{f}{\lambda_f}
+\frac{\lambda_g}{\lambda}\frac{g}{\lambda_g}\right) \\
&\le\left(\frac{\lambda_f}{\lambda}\right)^{p_-} m_{p(\cdot)}(f/\lambda_f)
\,+\,\left(\frac{\lambda_g}{\lambda}\right)^{p_-} m_{p(\cdot)}(g/\lambda_g) \\
&\le\left(\frac{\lambda_f}{\lambda}\right)^{p_-} + \left(\frac{\lambda_g}{\lambda}\right)^{p_-} 
\le 2^{1-p_-}\left(\frac{\lambda_f}{\lambda}+\frac{\lambda_g}{\lambda}\right)^{p_-} = 2^{1-p_-}.
\end{align*}
Then if follows by applying~\eqref{const-in-modular} with $\alpha=2^{1-1/p_-}$ that 
\[
m_{p(\cdot)}\left(\frac{f+g}{2^{1/p_- -1}\lambda}\right)\le2^{p_- -1}m_{p(\cdot)}\left(\frac{f+g}{\lambda}\right)\le1,
\] 
and thus $\|f+g\|_{p(\cdot)}\le2^{1/p_- -1}(\lambda_f+\lambda_g)$. Taking now infimum over all such $\lambda_f$
and $\lambda_g$, we get the desired quasi-triangle inequality with $C_\D=2^{1/p_- -1}$. Since, moreover, the lattice 
property for $\|\cdot\|_{p(\cdot)}$ holds as a consequence of the order preserving property of the modular,
we have proved that $L^{p(\cdot)}(\Omega,d,\mu)$ is a quasi-normed lattice with respect to $\|\cdot\|_{p(\cdot)}$. 

To complete the proof, let us verify the Fatou property. Take a sequence $\{f_n\}$ in $L^{p(\cdot)}(\Omega,\mu)$
such that $0\le f_n\uparrow f$ and 
\[
\Lambda:=\sup_{n\ge0}\|f_n\|_{p(\cdot)}<\infty.
\] 
With this choice, $m_{p(\cdot)}(f_n/\Lambda)\le1$ for any $n\ge0$. Using the
monotone convergence theorem~\cite[Theorem~2.14]{F99}, we get
\begin{align*}
m_{p(\cdot)}\left(\frac{f}{\Lambda}\right)&=
\int_{\Omega\backslash\Omega_\infty} \left|\frac{f(x)}{\Lambda}\right|^{p(x)}d\mu(x)+
\operatornamewithlimits{ess\,sup}_{x\in\Omega_\infty}\left|\frac{f(x)}{\Lambda}\right| \\
&=\lim_{n\to\infty} \left(\int_{\Omega\backslash\Omega_\infty} \left|\frac{f_n(x)}{\Lambda}\right|^{p(x)}d\mu(x)+
\operatornamewithlimits{ess\,sup}_{x\in\Omega_\infty}\left|\frac{f_n(x)}{\Lambda}\right|\right) \\
&=\lim_{n\to\infty}m_{p(\cdot)}\left(\frac{f_n}{\Lambda}\right)\le1,
\end{align*}
which implies that $f\in L^{p(\cdot)}(\Omega,\mu)$ and 
$\|f\|_{p(\cdot)}\le\sup_{n\ge0}\|f_n\|_{p(\cdot)}$. Thus, the Fatou property holds with the constant 
$C_\mathcal{F}=1$, and $L^{p(\cdot)}(\Omega,\mu)$ is a quasi-Banach lattice with the Luxemburg--Nakano quasi-norm. 
\end{proofof}

Hence, any space $L^{p(\cdot)}(\Omega,\mu)$ with $p_->0$ is a lattice and can be convexified; in fact, its every
convexification is a variable Lebesgue space itself---with the exponent function multiplied by the parameter of 
convexification. The following theorem extends \cite[Proposition~2.18]{CF13}.
\begin{theorem}\label{th:equiv-QN}
Given $s>0$ and an exponent function $p(\cdot)$ with $p_->0$, the $s$-convexification of $L^{p(\cdot)}(\Omega,\mu)$
coincides with $L^{sp(\cdot)}(\Omega,\mu)$, and $\|\cdot\|_{(L^{p(\cdot)})^{(s)}}$ and $\|\cdot\|_{sp(\cdot)}$ 
are equivalent quasi-norms, such that
\[
2^{-(1/s)\max\{1/p_-,1\}}\|f\|_{(L^{p(\cdot)})^{(s)}}\le\|f\|_{sp(\cdot)}\le
2^{\max\{1/sp_-,1\}}\|f\|_{(L^{p(\cdot)})^{(s)}}
\]
for all functions $f\in L^{sp(\cdot)}(\Omega,\mu)$.
\end{theorem}
\begin{proof}
Take an $f\in(L^{p(\cdot)})^{(s)}(\Omega,\mu)$. By the definition of the convexification, 
$|f|^s\in L^{p(\cdot)}(\Omega,\mu)$, and we fix an arbitrary $\lambda>\|\,|f|^s\|_{p(\cdot)}$. 
This choice of $\lambda$ implies that $m_{p(\cdot)}(|f|^s/\lambda)\le1$ and thus
\[
\int_{\Omega\setminus\Omega_\infty}\left(\frac{|f(x)|^s}{\lambda}\right)^{p(x)}d\mu(x)\le1,
\quad
\operatornamewithlimits{ess\,sup}_{x\in\Omega_\infty}\frac{|f(x)|^s}{\lambda}\le1,
\]
or equivalently,
\[
\int_{\Omega\setminus\Omega_\infty}\left(\frac{|f(x)|}{\lambda^{1/s}}\right)^{sp(x)}d\mu(x)\le1,
\quad
\operatornamewithlimits{ess\,sup}_{x\in\Omega_\infty}\frac{|f(x)|}{\lambda^{1/s}}\le1.
\]
Summing the last two inequalities gives $m_{sp(\cdot)}(f/\lambda^{1/s})\le2$, and so
$f\in L^{sp(\cdot)}(\Omega,\mu)$. It follows then by inequality~\eqref{const-in-modular} that
\[
m_{sp(\cdot)}\left(2^{-\max\{1/sp_-,1\}}\frac{f}{\lambda^{1/s}}\right)\le
\frac12 m_{sp(\cdot)}\left(\frac{f}{\lambda^{1/s}}\right)\le1,
\]
whence $\|f\|_{sp(\cdot)}\le2^{\max\{1/sp_-,1\}}\lambda^{1/s}$. Taking the infimum over all such $\lambda$ yields
\begin{equation}\label{ineq:half-1}
\|f\|_{sp(\cdot)}\le2^{\max\{1/sp_-,1\}}\|f\|_{(L^{p(\cdot)})^{(s)}}.
\end{equation}

For the converse, take a function $f\in L^{sp(\cdot)}(\Omega,\mu)$ and fix any $\lambda>\|f\|_{sp(\cdot)}$ 
noting that ${m_{sp(\cdot)}(f/\lambda)\le1}$. Arguing as we did above, we consequently obtain that 
$m_{p(\cdot)}(|f|^s/\lambda^s)\le2$. Further application of inequality~\eqref{const-in-modular} gives us
\[
m_{p(\cdot)}\left(2^{-\max\{1/p_-,1\}}\frac{|f|^s}{\lambda^s}\right)\le1,
\] 
which implies that $f\in(L^{p(\cdot)})^{(s)}(\Omega,\mu)$ and 
$\|\,|f|^s\|_{p(\cdot)}\le2^{\max\{1/p_-,1\}}\lambda^s$. By taking the infimum over all $\lambda$ considered, 
we obtain
\begin{equation}\label{ineq:half-2}
2^{-(1/s)\max\{1/p_-,1\}}\|f\|_{(L^{p(\cdot)})^{(s)}}\le\|f\|_{sp(\cdot)}.
\end{equation}
Inequalities~\eqref{ineq:half-1} and \eqref{ineq:half-2} are the two parts of the desired inequlity. 
\end{proof}

Let us remark, however, that the modular $m_{p(\cdot)}$ which we use to define the variable Lebesgue space 
$L^{p(\cdot)}(\Omega,\mu)$ is not the only functional suitable for this purpose. As an alternative, one may take the
``max''-modular
\[
m_{p(\cdot)}^\text{max}(f):=\max\left\{
\int_{\Omega\setminus\Omega_\infty}|f(x)|^{p(x)}d\mu(x),\ 
\operatornamewithlimits{ess\,sup}\limits_{x\in\Omega_\infty}|f(x)|
\right\} 
\]
and similarly declare $L^{p(\cdot)}(\Omega,\mu)$ to be the space of measurable functions $f$ for which 
$m_{p(\cdot)}^\text{max}(f/\lambda)<\infty$ at some $\lambda>0$. Clearly, the new definition covers the same
collection of functions as before. The analogue of the Luxemburg--Nakano quasi-norm based on the ``max''-modular, i.e.
\[
\|f\|_{p(\cdot)}^\text{max}:=\inf\{\lambda>0\ :\ m_{p(\cdot)}^\text{max}(f/\lambda)\le1\},
\]
is again a quasi-norm with the constant $C_\D=\max\{1,2^{1/p_- -1}\}$ whenever $p_->0$; in fact, the entire 
Theorem~\ref{th:Q-B-lattice} remains valid for $\|\cdot\|_{p(\cdot)}^\text{max}$---one can easily trace that
the proof goes identically line by line, with the corresponding changes from the sum to maximum where necessary, 
given that Lemma~\ref{le:elem-prop} still holds when $m_{p(\cdot)}$ is replaced by $m_{p(\cdot)}^\text{max}$.  

The two quasi-norms $\|\cdot\|_{p(\cdot)}$ and $\|\cdot\|_{p(\cdot)}^\text{max}$ on $L^{p(\cdot)}(\Omega,\mu)$
are equivalent: the simple inequality $\max\{a,b\}\le a+b\le2\max\{a,b\}$ provides the estimate
\[
\|f\|_{p(\cdot)}^\text{max}\le\|f\|_{p(\cdot)}\le2^{\max\{1/p_-,1\}}\|f\|_{p(\cdot)}^\text{max}.
\]
A curious observation about the maximum quasi-norm is that $\|\cdot\|_{p(\cdot)}^\text{max}$ is equal to 
the ``convexified'' quasi-norm based on it, while for the sum quasi-norm, we only had the weaker 
Theorem~\ref{th:equiv-QN}.
\begin{lemma}\label{le:equal-QN}
Given $s>0$ and an exponent function $p(\cdot)$ with $p_->0$, for any 
$f\in L^{sp(\cdot)}(\Omega,\mu)$ the \mbox{``convexified''} quasi-norm
\[
\|f\|_{(L^{p(\cdot)})^{(s)}}^\textnormal{max}:=(\|\,|f|^s\|_{p(\cdot)}^\textnormal{max})^{1/s}
\]
is equal to the Luxemburg--Nakano quasi-norm $\|f\|_{sp(\cdot)}^\textnormal{max}$.
\end{lemma}
\begin{proof}
It is easy to see that since for any $\lambda>0$,
\[
\operatornamewithlimits{ess\,sup}_{x\in\Omega_\infty}\frac{|f(x)|^s}{\lambda}\le1
\quad\text{if and only if}\quad
\operatornamewithlimits{ess\,sup}_{x\in\Omega_\infty}\frac{|f(x)|}{\lambda^{1/s}}\le1,
\]
there holds the set equality
\[
\left\{\lambda>0\ :\ m_{p(\cdot)}^\text{max}\left(\frac{|f|^s}{\lambda}\right)\le1\right\}=
\left\{\lambda>0\ :\ m_{sp(\cdot)}^\text{max}\left(\frac{f}{\lambda^{1/s}}\right)\le1\right\}.
\]
By passing to the infima of the two sets, we arrive at 
$\|\,|f|^s\|_{p(\cdot)}^\text{max}=(\|f\|_{sp(\cdot)}^\text{max})^s$, which yields the desired result. 
\end{proof}

In our experience, the maximum quasi-norm, though being less used than its ``sum'' counterpart, often 
brings more technical convenience to the proofs; so it is generally a good idea to check any 
$L^{p(\cdot)}$-related results for both.

\section{Review on Spaces of Homogeneous Type}\label{sec:SHT}

\subsection{Quasi-metric spaces}
First of all, spaces of homogeneous type are quasi-metric spaces, and thus have a topological structure 
weaker than metric spaces. By definition, a quasi-metric on a set $\Omega$ is a function 
$\rho:\Omega\times\Omega\to[0,\infty)$ such that 
\begin{enumerate}
\item[(1)] $\rho(x,y)=0$ if and only if $x=y$,

\item[(2)] $\rho(x,y)=\rho(y,x)$ for every $x,y\in\Omega$,

\item[(3)] there exists a constant $A_0>0$ such that for every $x,y,z\in\Omega$,
\[
\rho(x,y)\le A_0(\rho(x,z)+\rho(z,y)).
\] 
\end{enumerate} 
A pair $(\Omega,\rho)$ is called a quasi-metric space. Note that if the set $\Omega$ has at least two distinct points, then
necessarily the quasi-metric constant $A_0\ge1$. The case $A_0=1$, however, is the usual metric case.  

Once we have a notion of ``distance,'' we can define a ball 
\[
B_\rho(x,r):=\{y\in\Omega:\rho(x,y)<r\}
\]
centered at $x\in\Omega$ and of radius $r>0$. A quasi-metric $\rho$ naturally induces the topology $\tau_\rho$ on
$\Omega$, in which a set $G\subset\Omega$ is defined to be open if for each $x\in G$ there exists $\varepsilon>0$ 
such that $B_\rho(x,\varepsilon)\subset G$. 

An important fact about this canonical topology on quasi-metric spaces
is that the same topology $\tau_\rho$ is unequivocally induced by any other quasi-metric $\rho'$ equivalent to $\rho$ in the 
sence that there exists a constant $c\ge1$ such that for all $x,y\in\Omega$, 
\[
c^{-1}\rho(x,y)\le\rho'(x,y)\le c\rho(x,y).
\]
Indeed, for any $x\in\Omega$ and $\varepsilon>0$, such a relation implies that 
$B_{\rho'}(x,\varepsilon/c)\subset B_\rho(x,\varepsilon)$ and $B_\rho(x,\varepsilon/c)\subset B_{\rho'}(x,\varepsilon)$,
and these inclusions of balls yield $\tau_\rho=\tau_{\rho'}$.

Quasi-metric spaces are also known for a certain misfortune with the balls: when $A_0>1$,
some balls $B_\rho(x,r)$ may fail to be open. Hyt\"onen and Kairema gave an elegant and simple example
of such an occurrence in~\cite[p.~5]{HK12}: by considering $\Omega=\{-1\}\cup[0,\infty)$ with the usual distance
between all other pairs of points except $\rho(-1,0):=1/2$, one gets a quasi-metric space with $A_0\ge2$, in which the ball 
$B_\rho(-1,1)=\{-1,0\}$ does not contain any ball of the form $B_\rho(0,\varepsilon)$, and hence is not open. Another 
example can be found in~\cite[p.~4310]{PS09}.

Fortunately, Mac\'ias and Segovia \cite[Theorem~2]{MS79} proved that given any quasi-metric $\rho$, there exists an 
equivalent quasi-metric $d$ such that $d^\alpha$, for some $0<\alpha<1$, is a genuine metric. Since every ball $B_d(x,r)$
coincides then with a metric ball $B_{d^\alpha}(x,r^\alpha)$, it follows that $\tau_\rho=\tau_d=\tau_{d^\alpha}$ and
every $d$-ball is open in this common topology. From now on, we will assume that an arbitrary quasi-metric $\rho$ on 
$\Omega$ has been replaced with an equivalent well-behaved quasi-metric $d$ and benefit from the fact that the associated
balls $B(x,r):=B_d(x,r)$ are all open sets. 

\subsection{Setting for analysis}
Now that the ``ball-friendly'' quasi-metric is chosen, we define the setting of spaces of homogeneous type following 
Coifman and Weiss~\cite{CW71,CW77}, who introduced this notion in the 1970s. 
\begin{definition}
A space of homogeneous type $(\Omega,d,\mu)$ is a quasi-metric space $(\Omega,d)$ equipped with a nonnegative Borel 
measure $\mu$ satisfying the doubling condition 
\begin{equation}\label{eq:doubling}
0<\mu(B(x,r))\le A\mu(B(x,r/2))<\infty
\end{equation}\label{eq:doubling}
with an absolute constant $A:=A_\mu\ge1$ for all balls $B(x,r)$.
\end{definition}
By assuming that balls have positive, finite measure, we avoid trivial measures and ensure that $\mu$ is 
$\sigma$-finite.

The doubling condition~\eqref{eq:doubling} of a measure $\mu$ implies the following homogeneity property of 
the quasi-metric space $(\Omega,d)$: there is a natural number 
\[
A_1\le A^{3\log_2{A_0}+5}
\] 
such that any ball $B(x,r)$ contains at most $A_1$ points $x_i$ satisfying $d(x_i,x_j)\ge r/2$. This is the very first thing 
pointed out by Coifman and Weiss in their discussion of spaces of homogeneous type~\cite[p.~67]{CW71}; in fact, their 
original---more general---defintion of these spaces required the homogeneity condition instead of the existence of a 
doubling measure, which is reflected in the name of the spaces. 

The homogeneity condition, in turn, implies the geometric doubling property of the quasi-metric $d$ with the same 
constant $A_1$---namely, any ball $B(x,r)$ can be covered by at most $A_1$ balls of radius $r/2$ (since we may always 
choose the smaller balls to be centered at the evenly spread points $x_i$). Thus, every space of homogeneous type is
geometrically doubling; having said this, we gain access to an important tool of adjacent dyadic grids developed by
Hyt\"onen and Kairema for geometrically doubling quasi-metric spaces.

The usefulness of dyadic objects in harmonic analysis on the Euclidean space $\mathbb{R}^n$ has long been known,
and development of similar systems in a more general setting of quasi-metric spaces was only a matter of time. One of the 
first advancements in this area is due to Christ~\cite[Theorem~11]{C90a} (see also \cite[Ch.~VI, Theorem~14]{C90b}),
who constructed a system of sets on $(\Omega,d,\mu)$ satisfying many properties of the Euclidean dyadic cubes. 
His construction was further refined by Hyt\"onen and Kairema~\cite[Theorem~2.2]{HK12}; moreover, these authors 
observed that a number of results in harmonic analysis exploit several adjacent grids of dyadic cubes instead of just 
one fixed grid and designed an analogous system of dyadic grids in quasi-metric spaces with geometric 
doubling~\cite[Theorem~4.1]{HK12}. We will use the version of this result from~\cite[Theorem~4.1]{AHT17}.
\begin{theorem}\label{th:Hytonen-Kairema}
Let $(\Omega,d)$ be a quasi-metric space with the quasi-metric constant $A_0$ and satisfy the geometric doubling 
condition with the constant $A_1$. Suppose the parameter $\delta\in(0,1)$ satisfies $96A_0^6\delta\le 1$. Then there 
exist countable sets of points $\{z_\alpha^{k,t}:\alpha\in\mathcal{A}_k\}$, $k\in\mathbb{Z}$, 
$t=1,2,\dots,K=K(A_0,A_1,\delta)$, and a finite number of dyadic grids 
$\mathcal{D}^t:=\{Q_\alpha^{k,t}:k\in\mathbb{Z},\alpha\in\mathcal{A}_k\}$, such that
\begin{enumerate}
\item[{\rm(a)}]
for every $t\in\{1,2,\dots,K\}$ and $k\in\mathbb{Z}$ one has
\begin{enumerate}
\item[{\rm (i)}]
$\Omega=\bigcup_{\alpha\in\mathcal{A}_k} Q_\alpha^{k,t}$ (disjoint union);

\item[{\rm (ii)}]
if $Q,P\in\mathcal{D}^t$, then $Q\cap P\in\{\emptyset, Q,P\}$;

\item[{\rm (iii)}]
if $Q_\alpha^{k,t}\in\mathcal{D}^t$, then
\begin{equation}\label{eq:dyadic-cubes}
B(z_\alpha^{k,t},c_1\delta^k)
\subset 
Q_\alpha^{k,t}
\subset 
B(z_\alpha^{k,t},C_1\delta^k),
\end{equation}
where $c_1:=(12A_0^4)^{-1}$ and $C_1:=4A_0^2$;

\item[{\rm(iv)}]
if $Q_\alpha^{k,t}\in\mathcal{D}^t$, there exist at least one $Q_\beta^{k+1,t}\in\mathcal{D}^t$, 
which is called a child of $Q_\alpha^{k,t}$, and exactly one $Q_\gamma^{k-1,t}\in\mathcal{D}^t$, 
which is called the parent of $Q_\alpha^{k,t}$, such that 
$Q_\beta^{k+1,t}\subset Q_\alpha^{k,t}\subset Q_\gamma^{k-1,t}$;
\end{enumerate}

\item[{\rm (b)}]
for every ball $B=B(x,r)$, there exists a cube $Q_B\in\bigcup_{t=1}^K\mathcal{D}^t$ such that $B\subset Q_B$ and
 $Q_B=Q_\alpha^{k-1,t}$ for some indices $\alpha\in\mathcal{A}_k$ and $t\in\{1,\dots,K\}$, where $k$ is 
the unique integer satisfying $\delta^{k+1}<r\le\delta^k$.
\end{enumerate}
\end{theorem}
For further work, fix a collection of the adjacent dyadic grids $\{\mathcal{D}^t:t=1,2,\ldots,K\}$ defined in 
Theorem~\ref{th:Hytonen-Kairema} and denote their union by
\[
\mathcal{D}:=\bigcup_{t=1}^K\mathcal{D}^t.
\]
The sets $Q^{k,t}_\alpha\in\mathcal{D}$ are Borel sets~\cite[Remark~4.2]{AHT17} referred to as dyadic cubes with 
centers $z^{k,t}_\alpha$ and sidelengths $\ell(Q^{k,t}_\alpha)=\delta^k$. Of course, these are not cubes in the 
standard Euclidean sense---but the name was preserved given the similar properties of each individual dyadic system 
$\mathcal{D}^t$, listed in (a), to those of the classical dyadic grid in $\mathbb{R}^n$. In particular,
inclusion~\eqref{eq:dyadic-cubes} implies that for any cube $Q^{k,t}_\alpha$, there are the containing ball
$B(z_\alpha^{k,t},C_1\delta^k)$ and the contained ball $B(z_\alpha^{k,t},c_1\delta^k)$, so we can trap a cube
``between'' two balls just like in $\mathbb{R}^n$. Conversely, part (b) guarantees that any quasi-metric ball is 
contained in a dyadic cube of a comparable size. 

More specifically, in spaces of homogeneous type---where we have a doubling measure $\mu$---dyadic
properties~(iii) and (b) yield the following relations between measures of cubes and their containing balls, and vice 
versa~\cite[Corollary~7.4]{HK12}.
\begin{lemma}\label{le:balls-cubes}
There exists a constant $C=C(A_0,\delta)\ge1$ such that for every dyadic cube $Q\in\mathcal{D}$ we have 
$\mu(B_Q)\le C\mu(Q)$, where $B_Q$ is the containing ball of $Q$ as in~\eqref{eq:dyadic-cubes}. Conversely, 
given a ball $B=B(x,r)$, there exist a dyadic grid $\mathcal{D}^t$ and a dyadic cube $Q_B\in\mathcal{D}^t$ 
such that $B\subset Q_B$ and $\mu(Q_B)\le C\mu(B)$.
\end{lemma}
This result is crucial for establishing equivalence between the classical maximal operator and its ``dyadic'' 
counterpart. 

\subsection{Maximal operators}\label{subsec:MOs} 
Let $f$ be a measurable function on $(\Omega,d,\mu)$, from now on denoted by 
$f\in L^0(\Omega,d,\mu)$. 
Then $Mf$, the Hardy--Littlewood maximal function of $f$, is defined for any
$x\in\Omega$ by
\[
Mf(x):=\sup_{B\ni x}\frac1{\mu(B)}\int_B |f(y)|d\mu(y),
\]
where the supremum is taken over all quasi-metric balls $B\subset\Omega$ containing $x$. The maximal operator
$M$ is a sublinear operator acting by the rule $f\mapsto Mf$. Simply put, the function $Mf$ renders the ``maximum'' 
mean value of $f$ about the point $x$. By maximizing the means of other orders $r>0$, we also define the 
parameterized maximal operator $M_r$ as
\begin{equation}\label{eq:M_r}
M_r f(x):=M(|f|^r)(x)^{1/r}=\sup_{B\ni x}\left(\frac1{\mu(B)}\int_B |f(y)|^r d\mu(y)\right)^{1/r}.
\end{equation}
Obviously, this formula yields the usual maximal operator $M$ when $r=1$.

As the following lemma shows, values of the parameterized maximal function increase with the parameter. 
\begin{lemma}\label{le:Mr-order}
If $0<r<s$ and $f\in L^0(\Omega,d,\mu)$, then for every $x\in\Omega$ there holds
\[
M_r f(x)\le M_s f(x).
\]
\end{lemma}
\begin{proof}
Fix an $x\in\Omega$ and a ball $B$ containing $x$. By H\"older's inequality with the exponent $s/r>1$, we have
\begin{align*}
\frac1{\mu(B)}\int_B |f(y)|^r d\mu(y) &=
\frac1{\mu(B)}\int_\Omega |f(y)|^r \chi_B(y) d\mu(y) \\
&\le \frac1{\mu(B)} \left(\int_B |f(y)|^s d\mu(y)\right)^{r/s} \mu(B)^{1-r/s} \\
&= \left(\frac1{\mu(B)} \int_B |f(y)|^s d\mu(y)\right)^{r/s} \le M_s f(x)^r.
\end{align*}
This implies that 
\[
M_r f(x)=\sup_{B\ni x} \left(\frac1{\mu(B)}\int_B |f(y)|^r d\mu(y)\right)^{1/r}\le M_s f(x).
\]
\end{proof}

Along with the classical maximal operator $M$, we will use the ``dyadic'' maximal operator $M^\mathcal{D}$ defined
for $f\in L^0(\Omega,d,\mu)$ by 
\[
M^\mathcal{D}f(x)=\sup_{Q\ni x} \frac1{\mu(Q)}\int_Q |f(y)|d\mu(y),\quad x\in\Omega,
\]
where the supremum is taken over all dyadic cubes $Q\in\mathcal{D}$ containing $x$. In style of Anderson, Hyt\"onen and
Tapiola~\cite{AHT17}, we put quotation marks around the word dyadic when we want to emphasize that the definition
in question uses the collection $\mathcal{D}$ of adjacent dyadic grids instead of a single grid. 

Note that we can similarly define the ``dyadic'' parameterized maximal operator $M_r^\mathcal{D}$ by taking the
supremum over cubes instead of balls in~\eqref{eq:M_r}. Lemma~\ref{le:Mr-order} still remains true if we replace 
$M_r$ and $M_s$ by their ``dyadic'' versions $M_r^\mathcal{D}$ and $M_s^\mathcal{D}$.

Due to Lemma~\ref{le:balls-cubes}, we easily deduce the pointwise equivalence of the maximal functions $Mf$ and 
$M^\mathcal{D}f$ (cf.~\cite[Proposition~7.9]{HK12}).
\begin{theorem}\label{th:M-dyadicM}
Let $f\in L^0(\Omega,d,\mu)$. For all $x\in\Omega$, we have the poinwise estimates 
\[
M^{\mathcal{D}}f(x)\le CMf(x)
\quad\text{and}\quad
Mf(x)\le CM^{\mathcal{D}}f(x)
\]
with the same constant $C\ge1$, independent of $f$, as in Lemma~\ref{le:balls-cubes}. 
\end{theorem}
\begin{proof}
Assume that $x\in Q$, $Q\in\mathcal{D}$. Let $B$ be the containing ball of $Q$ as in~\eqref{eq:dyadic-cubes}. Then 
$\mu(B)\le C\mu(Q)$ by Lemma~\ref{le:balls-cubes}, and therefore
\[
\frac1{\mu(Q)}\int_Q |f(y)|d\mu(y) \le \frac{C}{\mu(B)}\int_B |f(y)|d\mu(y) \le CMf(x).
\]
By taking the supremum over all dyadic cubes in $\mathcal{D}$ containing $x$, we conclude from the above inequality that
$M^{\mathcal{D}}f(x)\le CMf(x)$.

For the reverse inequality, consider a ball $B\ni x$ and let $Q\in\mathcal{D}^t$, for some $t$, be the dyadic cube from 
Lemma~\ref{le:balls-cubes} such that $B\subset Q$ and $\mu(Q)\le C\mu(B)$. Then, just like in the argument above, 
but with the interchanged roles of $B$ and $Q$, we have
\[
\frac1{\mu(B)}\int_B |f(y)|d\mu(y) \le \frac{C}{\mu(Q)}\int_Q |f(y)|d\mu(y) \le CM^{\mathcal{D}}f(x),
\]
which results in $Mf(x)\le CM^{\mathcal{D}}f(x)$ after taking the supremum over all balls $B$ cointaining $x$.
\end{proof}

Both ``dyadic'' and classical maximal operators are countably subadditive, as stated in the next lemma.
\begin{lemma}\label{le:count-subadd}
If $\{f_k\}\subset L^0(\Omega,d,\mu)$ and $f=\sum_{k=0}^\infty f_k$, then for every $x\in\Omega$,
\[
M^\mathcal{D} f(x)\le\sum_{k=0}^\infty M^\mathcal{D}f_k(x).
\] 
The same is true if we replace $M^\mathcal{D}$ by $M$.
\end{lemma}
\begin{proof}
The statement is a direct consequence of the countable additivity of the integral (see~\cite[Theorem~2.15]{F99}) 
resulting from the monotone convergence theorem---since for any fixed $x\in\Omega$, we have
\begin{align*}
M^\mathcal{D}f(x)&=\sup_{Q\ni x} \frac1{\mu(Q)} \int_Q |f(y)|d\mu(y)
\le \sup_{Q\ni x} \frac1{\mu(Q)} \int_Q \sum_{k=0}^\infty |f_k(y)|d\mu(y) \\
&=\sup_{Q\ni x} \sum_{k=0}^\infty \frac1{\mu(Q)} \int_Q |f_k(y)|d\mu(y) 
\le \sum_{k=0}^\infty M^\mathcal{D} f_k(x),
\end{align*}
where all the suprema are taken over cubes $Q\in\mathcal{D}$ containing $x$. When replacing dyadic cubes by balls 
$B\ni x$, we obtain the same inequality for the operator $M$.
\end{proof}

\subsection{``Dyadic'' $A_1$ weights}\label{subsec:A1-weights}
A nonnegative measurable function $w$ on $\Omega$ is said to be a weight. The ``dyadic'' class $A_1^\mathcal{D}$
consists of the weights $w$ such that 
\[
[w]_{A_1^\mathcal{D}}:=\operatornamewithlimits{ess\,sup}\limits_{x\in\Omega} 
\frac{M^\mathcal{D}w(x)}{w(x)}<\infty.
\] 

For each $Q\in\mathcal{D}$, define the set of all its nearby cubes of that same generation, 
\[
\mathcal{N}_Q:=\{Q'\in\mathcal{D}\ :\ Q'\cap Q\ne\emptyset,\;\ell(Q')=\ell(Q)\}.
\]
According to \cite[Definition~4.4]{AHT17}, a generalized dyadic parent (gdp) of $Q$ is any cube $Q^*$ such that 
$\ell(Q^*)=\delta^{-2}\ell(Q)$ and for every $Q'\in\mathcal{N}_Q$ we have $Q'\subset Q^*$. Every cube has 
at least one gdp (see~\cite[Lemma~4.5]{AHT17}), and we set $S\ge1$ to be such a constant that for all 
$Q\in\mathcal{D}$ and $Q'\in\mathcal{N}_Q$, there holds $\mu(Q^*)\le S\mu(Q')$. This constant, introduced 
in~\cite[Section~5]{AHT17}, is always finite---which we now prove by giving one of its possible values.
\begin{theorem}
Given that $S=A(A_0/\delta^3)^{\log_2 A}$, there holds $\mu(Q^*)\le S\mu(Q')$ for any gdp $Q^*$ and 
all $Q'\in\mathcal{N}_Q$ related to an arbitrary cube $Q\in\mathcal{D}$.  
\end{theorem}
\begin{proof}
Note first that the doubling property~\eqref{eq:doubling} of the measure $\mu$ implies that for all $x\in\Omega$, 
$0<r<R$ and $y\in B(x,R)$, one has
\begin{equation}\label{eq:lower-mass-bound}
\frac{\mu(B(x,R))}{\mu(B(y,r))}\le A\left(\frac{2A_0 R}{r}\right)^{\log_2 A}
\end{equation}
(cf.~\cite[Lemma~2.3]{CS18}). This follows by showing that for all $z\in B(x,R)$, the quasi-triangle inequality yields 
\[
d(y,z)\le A_0(d(y,x)+d(x,z))<2A_0 R
\]	
and consequently $B(x,R)\subset B(y,2A_0 R)$, and then applying the well-known relation between the measures of 
concentric balls, see~\cite[Lemma~7.3]{HK12}, to counclude that
\[
\mu(B(x,R))\le\mu(B(y,2A_0 R))\le A\left(\frac{2A_0 R}{r}\right)^{\log_2 A}\mu(B(y,r)).
\]  

Let $\ell(Q)=\ell(Q')=\delta^k$; then $\ell(Q^*)=\delta^{k-2}$. According to Property~(a)-(iii) from 
Theorem~\ref{th:Hytonen-Kairema}, there exist $x,y\in\Omega$ such that 
\[
B(y,c_1\delta^k)\subset Q'\subset Q^*\subset B(x,C_1 \delta^{k-2}).
\] 
If we now apply inequality~\eqref{eq:lower-mass-bound}, recall the values of $c_1$ and $C_1$, and make use of 
the definition of the parameter $\delta$ satisfying $96A_0^6\delta\le1$, we get that
\begin{align*}
\mu(Q^*)&\le\mu(B(x,C_1\delta^{k-2})) \\
&\le A\left(\frac{2A_0 C_1 \delta^{k-2}}{c_1\delta^k}\right)^{\log_2 A} \mu(B(y,c_1\delta^k)) \\
&\le A\left(\frac{A_0}{\delta^3}\right)^{\log_2 A}\mu(Q'),
\end{align*}
and thus the claim follows.
\end{proof}

The following result is an easy consequence of the weak reverse H\"older inequality from~\cite[Theorem~5.4]{AHT17}.
We present it as given in the work by Karlovich~\cite[Lemma~6]{K19}, and also take this occasion to correct a misleading 
typo in~\cite{K19}: throughout Section~3 there, $\mathcal{D}$ is meant to be the union of adjacent dyadic grids, just like
in our text, instead of a single dyadic grid.
\begin{theorem}
Let $w\in A_1^{\mathcal{D}}$ and $Q\in\mathcal{D}$. Then for every $\eta$ satisfying
\begin{equation}\label{eq:eta-range}
0<\eta\le\frac1{2S^2 K[w]_{A_1^\mathcal{D}}},
\end{equation}
one has 
\begin{equation}\label{eq:RHI-2}
\left(\frac{1}{2\mu(Q)}\int_Q w^{1+\eta}(y)\,d\mu(y)\right)^{\frac{1}{1+\eta}}
\le 
S[w]_{A_1^\mathcal{D}}\frac{1}{\mu(Q)}\int_Q w(y)\,d\mu(y).
\end{equation}
\end{theorem}
As a corollary to this theorem, we have a pointwise estimate for the ``dyadic'' parameterized maximal functions of
$A_1^\mathcal{D}$ weights.  
\begin{corollary}\label{cor:point-est}
If $w\in A_1^\mathcal{D}$, then for every $\eta$ in the range~\eqref{eq:eta-range} and for a.\,e. $x\in\Omega$,
\begin{equation*}
M_{1+\eta}^\mathcal{D}w(x)\le2S[w]_{A_1^\mathcal{D}}^2 w(x).
\end{equation*}
\end{corollary}
\begin{proof}
Take an $x\in\Omega$ and fix a cube $Q\ni x$. Since $w\in A_1^\mathcal{D}$, it follows from inequality~\eqref{eq:RHI-2}
and the definition of the $A_1^\mathcal{D}$ weight that for all $\eta$ satisfying~\eqref{eq:eta-range},
\begin{align*}
\left(\frac{1}{\mu(Q)}\int_Q w^{1+\eta}(y)\,d\mu(y)\right)^{\frac{1}{1+\eta}}&\le
2^{\frac1{1+\eta}}S[w]_{A_1^\mathcal{D}}\frac1{\mu(Q)}\int_Q w(y)\,d\mu(y) \\
&\le 2S[w]_{A_1^\mathcal{D}} M^\mathcal{D}w(x) \\
&\le 2S[w]_{A_1^\mathcal{D}}^2 w(x),
\end{align*}
where the last inequality holds for almost every $x$. Then, if we take the supremum over all cubes $Q$ containing $x$ on
the left-hand side, we get the desired inequality.
\end{proof}

\section{Proof of the Main Theorem}\label{sec:proof}

\subsection{Boundedness on ``higher'' convexifications}
Consider first a quasi-metric measure space $(\Omega,d,\mu)$ without the requirement on $\mu$ 
to be doubling. We easily prove that the boundedness of $M$ on a quasi-normed lattice
over $(\Omega,d,\mu)$ entails its boundedness on all of the ``higher'' convexifications of the lattice.
\begin{lemma}\label{le:trivial-impl}
Let $X(\Omega,d,\mu)$ be a quasi-normed lattice over a quasi-metric measure space $(\Omega,d,\mu)$.
If $M$ is bounded on $X^{(r)}(\Omega,d,\mu)$ for some $r>0$, then $M$ is bounded on 
$X^{(s)}(\Omega,d,\mu)$ for all $s>r$ and
\begin{equation}\label{ineq:M-norm}
\|M\|_{X^{(s)}\to X^{(s)}}\le \|M\|_{X^{(r)}\to X^{(r)}}^{r/s}.
\end{equation}
\end{lemma}
\begin{proof}
Fix an $s>r$ and take an $f\in X^{(s)}(\Omega,d,\mu)$. Note that since $X(\Omega,d,\mu)$ is a quasi-normed lattice, 
so are $X^{(r)}(\Omega,d,\mu)$ and $X^{(s)}(\Omega,d,\mu)$ by Lemma~\ref{le:convex-lattice}. Then the pointwise
 inequality from Lemma~\ref{le:Mr-order} and the boundedness of $M$ on $X^{(r)}(\Omega,d,\mu)$ imply that
\begin{align*}
\|Mf\|_{X^{(s)}} &= \|(Mf)^s\|_X^{1/s} = \|(Mf)^{s/r}\|_{X^{(r)}}^{r/s} = \|M_{r/s}(|f|^{s/r})\|_{X^{(r)}}^{r/s} \\
&\le \|M(|f|^{s/r})\|_{X^{(r)}}^{r/s} \le \|M\|_{X^{(r)}\to X^{(r)}}^{r/s} \|\,|f|^{s/r}\|_{X^{(r)}}^{r/s} \\
&= \|M\|_{X^{(r)}\to X^{(r)}}^{r/s} \|f\|_{X^{(s)}}.
\end{align*}
Thus, $M$ is bounded on $X^{(s)}(\Omega,d,\mu)$ and \eqref{ineq:M-norm} is satisfied. 
\end{proof}

\subsection{Self-improving property}
Extending the boundedness of $M$ from a quasi-normed lattice to a range of its ``lower'' convexifications, however, 
is not at all trivial. We were able to construct such an extension for a more specific class of lattices---quasi-Banach
lattices with the Fatou property. In contrast to Lemma~\ref{le:trivial-impl}, the assumption that the lattices are taken 
over spaces of homogeneous type is crucial for the converse result.

As we noted in the Introduction, it is this converse result that is generally called the \textit{self-improving property} 
of the boundedness of the maximal operator. We establish it in two parts---first, ``pushing'' the parameter of 
convexification from $r=1$ down (Theorem~\ref{th:main-lower}), and then adapting the designed argument to pass 
from convexifications of order $s>1$ right down to $s=1$ (Theorem~\ref{th:main-upper}).
 
\begin{theorem}\label{th:main-lower}
Let $X(\Omega,d,\mu)$ be a quasi-Banach lattice with the Fatou property. Suppose $M$ is bounded on
$X(\Omega,d,\mu)$. Then there exists $r_0\in(0,1)$ such that if $r\in[r_0,1)$, then $M$ is bounded on 
$X^{(r)}(\Omega,d,\mu)$.
\end{theorem}
\begin{proof}
Instead of proving the theorem directly, we show that the same statement is true for the ``dyadic'' maximal operator 
$M^\mathcal{D}$ and use the equivalence of the maximal functions $Mf$ and $M^\mathcal{D}f$ to pass from $M$ to 
$M^\mathcal{D}$ and backwards. 

Since $M$ is bounded on $X(\Omega,d,\mu)$, it follows from Theorem~\ref{th:M-dyadicM} that $M^\mathcal{D}$ is also
bounded on $X(\Omega,d,\mu)$ and $\|M^\mathcal{D}\|_{X\to X}\le C\|M\|_{X\to X}$. To show that there is a number
$r_0\in(0,1)$ such that $M^\mathcal{D}$ is then bounded on any $X^{(r)}(\Omega,d,\mu)$ with $r_0\le r<1$, we first 
construct an $A_1^\mathcal{D}$ weight applying the Rubio de Francia iteration algorithm (see, e.g., \cite{CMP11} for a 
nice introduction to techniques based on it).

Fix an $\varepsilon$ such that $0<\varepsilon\|M^\mathcal{D}\|_{X\to X}<2^{-1/\rho}$, where $\rho$ is the 
Aoki-Rolewicz exponent of the lattice $X(\Omega,d,\mu)$. Given $h\in X(\Omega,d,\mu)$, define
\[
\mathcal{R}_\varepsilon^\mathcal{D}h(x):=\sum_{k=0}^\infty \varepsilon^k (M^\mathcal{D})^k h(x),
\quad x\in\Omega,
\]
where $(M^\mathcal{D})^k$ is the operator $M^\mathcal{D}$ iterated $k$ times and $(M^\mathcal{D})^0 h=|h|$. This
``dyadic'' Rubio de Francia operator has the following properties:
\begin{enumerate}
\item[(a)] trivially, $|h(x)|\le\mathcal{R}_\varepsilon^\mathcal{D}h(x)$ for all $x\in\Omega$;

\item[(b)] $\mathcal{R}_\varepsilon^\mathcal{D}$ is bounded on $X(\Omega,d,\mu)$ and 
$\|\mathcal{R}_\varepsilon^\mathcal{D}\|_{X\to X}\le 4^{1/\rho}C_\mathcal{F}$
(like before, $C_\mathcal{F}$ denotes the constant in the Fatou property for the lattice $X(\Omega,d,\mu)$);

\item[(c)] $\mathcal{R}_\varepsilon^\mathcal{D}h\in A_1^\mathcal{D}$ and 
$[\mathcal{R}_\varepsilon^\mathcal{D}h]_{A_1^\mathcal{D}}\le1/\varepsilon$.
\end{enumerate}

Property~(b) follows from Theorem~\ref{th:Aoki-Rol-inf}, according to which
\begin{align*}
\|\mathcal{R}_\varepsilon^\mathcal{D}h\|_X &\le 
2^{1/\rho}C_\mathcal{F} \left(\sum_{k=0}^\infty \varepsilon^{k\rho}\|(M^\mathcal{D})^k h\|_X^\rho\right)^{1/\rho}\\
&\le 2^{1/\rho}C_\mathcal{F} 
\left(\sum_{k=0}^\infty \left(\varepsilon\|M^\mathcal{D}\|_{X\to X}\right)^{k\rho}\right)^{1/\rho} \|h\|_X \\
&\le 2^{1/\rho}C_\mathcal{F} \left(\sum_{k=0}^\infty (2^{-1/\rho})^{k\rho}\right)^{1/\rho} \|h\|_X \\
&= 4^{1/\rho}C_\mathcal{F}\|h\|_X.
\end{align*}
Property~(c) is a consequence of the ``countable subadditivity'' (see Lemma~\ref{le:count-subadd}) and homogeneity 
of the maximal operator, since for any $x\in\Omega$ we have
\[
M^\mathcal{D}(\mathcal{R}_\varepsilon^\mathcal{D}h)(x) \le
\sum_{k=0}^\infty\varepsilon^k (M^\mathcal{D})^{k+1}h(x) \le
\frac1\varepsilon \mathcal{R}_\varepsilon^\mathcal{D}h(x).
\]

Denote $\eta_0:=\varepsilon/(2S^2K)$. By Property~(c), 
\[
\eta_0\le\frac{1}{2S^2K[\mathcal{R}_\varepsilon^\mathcal{D}h]_{A_1^\mathcal{D}}},
\]
therefore, it follows from Corollary~\ref{cor:point-est} and again~(c) that for every $0<\eta\le\eta_0$ and almost every
$x\in\Omega$, there holds 
\begin{equation}\label{eq:param-M-weight}
M_{1+\eta}^\mathcal{D} (\mathcal{R}_\varepsilon^\mathcal{D}h)(x)\le
\frac{2S}{\varepsilon^2}\mathcal{R}_\varepsilon^\mathcal{D}h(x).
\end{equation} 

Declare $r_0=1/(1+\eta_0)$. Then fix $r\in[r_0,1)$ and determine the unique number $\eta\in(0,\eta_0]$ such that 
$r=1/(1+\eta)$. Then by Properties~(a) and (b) and inequality~\eqref{eq:param-M-weight}, for any 
$f\in X^{(r)}(\Omega,d,\mu)$ we have
\begin{align*}
\|M^\mathcal{D}f\|_{X^{(r)}}&=\|(M^\mathcal{D}f)^r\|_X^{1/r}=\|M_{1/r}^\mathcal{D}(|f|^r)\|_X^{1/r} \\
&= \|M_{1+\eta}^\mathcal{D}(|f|^r)\|_X^{1/r} \le 
\|M_{1+\eta}^\mathcal{D}(\mathcal{R}_\varepsilon^\mathcal{D}(|f|^r))\|_X^{1/r} \\
&\le \left(\frac{2S}{\varepsilon^2}\right)^{1/r}
\|\mathcal{R}_\varepsilon^\mathcal{D}(|f|^r)\|_X^{1/r} \le
\left(\frac{2S}{\varepsilon^2}\right)^{1/r} (4^{1/\rho}C_\mathcal{F})^{1/r}\|\,|f|^r\|_X^{1/r} \\
&= \left(\frac{2S\cdot4^{1/\rho}C_\mathcal{F}}{\varepsilon^2}\right)^{1/r}\|f\|_{X^{(r)}},
\end{align*}
which means that $M^\mathcal{D}$ is bounded on $X^{(r)}(\Omega,d,\mu)$. Once again, Theorem~\ref{th:M-dyadicM}
guarantees that in this case $M$ is also bounded on all $X^{(r)}(\Omega,d,\mu)$, $r_0\le r<1$, and
$\|M\|_{X^{(r)}\to X^{(r)}}\le C\|M^\mathcal{D}\|_{X^{(r)}\to X^{(r)}}$.
\end{proof}
\begin{remark}
The construction in the proof implies that 
\[
r_0>\left(1+\frac{C}{2^{1/\rho+1}S^2 K\|M\|_{X\to X}}\right)^{-1}.
\]
Intuitively, this lower estimate for the parameter of convexification $r_0$ down to which the boundedness of $M$ can be 
extended---depending on the quasi-norm of the maximal operator---admits the following interpretation: a smaller
value of $\|M\|_{X\to X}$ allows to achieve a lower (that is, better) self-improvement ``threshold'' $r_0$.
\end{remark}

The second part of the self-improving property uses the argument from the above proof and follows almost as its 
corollary. 

\begin{theorem}\label{th:main-upper}
Let $X(\Omega,d,\mu)$ be a quasi-Banach lattice with the Fatou property. Suppose that there exists $s_0>1$ such that
$M$ is bounded on $X^{(s)}(\Omega,d,\mu)$ for every $1<s<s_0$, and 
\begin{equation}\label{eq:limit}
\lim_{s\to1^+} (s-1)\|M\|_{X^{(s)}\to X^{(s)}}=0.
\end{equation}
Then $M$ is bounded on $X(\Omega,d,\mu)$.
\end{theorem}
\begin{proof}
It is immediate from the equivalence of the classical and ``dyadic'' maximal functions, given by 
Theorem~\ref{th:M-dyadicM}, that not only $M$ but also $M^\mathcal{D}$ is bounded on each $X^{(s)}(\Omega,d,\mu)$,
$1<s<s_0$, and limit~\eqref{eq:limit} holds true if we replace $M$ by $M^\mathcal{D}$. 

Therefore, we can find an $s\in(1,s_0)$ such that 
\[
4C_\D(X)\cdot2S^2K(s-1)\|M^\mathcal{D}\|_{X^{(s)}\to X^{(s)}}<1
\]
and fix $\varepsilon=2S^2K(s-1)$. Denote by $\rho_s$ the Aoki-Rolewicz exponent of $X^{(s)}(\Omega,d,\mu)$; 
note that $2^{1/\rho_s}=2C_\D(X^{(s)})\le4C_\D(X)$ due to inequality~\eqref{eq:C-tri-convex}. Then 
\[
0<\varepsilon\|M^\mathcal{D}\|_{X^{(s)}\to X^{(s)}}<\frac1{4C_\D(X)}\le2^{-1/\rho_s}.
\]

After this choice of $\varepsilon$, we can repeat the argument from the proof of Theorem~\ref{th:main-lower} almost
verbatim---it suffices to apply it to $X=X^{(s)}$ and $r=r_0=1/s$ with the corresponding change of constants 
$\rho=\rho_s$ and $C_\mathcal{F}=C_\mathcal{F}(X^{(s)})$. Following this ``re-designation,'' $X^{(r)}$ becomes 
$(X^{(s)})^{(1/s)}=X^{(1)}=X$. By Lemma~\ref{le:Fatou}, $C_\mathcal{F}(X^{(s)})=C_\mathcal{F}(X)^{1/s}$,
and we conclude from the earlier argument that for all $f\in X(\Omega,d,\mu)$,
\[
\|M^\mathcal{D}f\|_X\le
C_\mathcal{F}(X) \left(\frac{2S\cdot4^{1/\rho_s}}{\varepsilon^2}\right)^s \|f\|_X,
\]
hence the ``dyadic'' operator $M^\mathcal{D}$ is bounded on $X(\Omega,d,\mu)$. By equivalence, $M$ is bounded 
on $X(\Omega,d,\mu)$, too.
\end{proof}

\subsection{Proof of Theorem~\ref{main-theorem}}
By bringing together the results of Lemma~\ref{le:trivial-impl} and Theorems~\ref{th:main-lower} and \ref{th:main-upper},
we obtain the main Theorem~\ref{main-theorem}. 

Indeed, implications $(1)\Rightarrow(2)$ and $(3)\Rightarrow(1)$
essentially give the extention of the boundedness of $M$ to the ``higher'' convexifications of a lattice, and are thus ensured
by Lemma~\ref{le:trivial-impl}. Implication $(1)\Rightarrow(3)$ is exactly Theorem~\ref{th:main-lower}. Finally, 
implication $(2)\Rightarrow(1)$ is a corollary of Theorem~\ref{th:main-upper}, since Condition~(2) taken as a premise can 
be weakened by assuming that there exists $s_0>1$ such that $M$ is bounded on $X^{(s)}(\Omega,d,\mu)$ if $1<s<s_0$
and limit~\eqref{eq:limit} holds. \qed

\subsection{Application to variable Lebesgue spaces}\label{subsec:appl}
Finally, we give an application of Theorem~\ref{main-theorem} to the spaces $L^{p(\cdot)}(\Omega,d,\mu)$ with 
$p_->0$, which by Theorem~\ref{th:Q-B-lattice} are quasi-Banach lattices with respect to the quasi-norm
$\|\cdot\|_{p(\cdot)}$ and have the Fatou property. To translate the 
applied result from the abstract language of convexifications to the clear language of variable Lebesgue spaces with 
different exponent functions, we use Theorem~\ref{th:equiv-QN} and therefore obtain the following special-case version 
of the main theorem. 
\begin{corollary}\label{cor:L^p}
Given a variable Lebesgue space $L^{p(\cdot)}(\Omega,d,\mu)$ with $p_->0$, the following are equivalent:
\begin{enumerate}
\item[\rm (1)] $M$ is bounded on $L^{p(\cdot)}(\Omega,d,\mu)$.

\item[\rm (2)] For all $s>1$, $M$ is bounded on $L^{sp(\cdot)}(\Omega,d,\mu)$ and 
\[
\lim_{s\to1^+}(s-1)\|M\|_{L^{sp(\cdot)}\to L^{sp(\cdot)}}=0.
\] 

\item[\rm (3)] There exists $r_0\in(0,1)$ such that if $r\in[r_0,1)$, then $M$ is bounded on $L^{rp(\cdot)}(\Omega,d,\mu)$.
\end{enumerate}
\end{corollary}
\begin{proof}
Once we set $X(\Omega,d,\mu)=L^{p(\cdot)}(\Omega,d,\mu)$ in the statement of Theorem~\ref{main-theorem}, the 
desired result follows immediately if the boundedness of $M$ on convexifications
$(L^{p(\cdot)})^{(s)}(\Omega,d,\mu)=L^{sp(\cdot)}(\Omega,d,\mu)$, $s>0$,
is understood with respect to the ``convexified'' operator quasi-norm defined by
\[
\|M\|_{(L^{p(\cdot)})^{(s)}\to(L^{p(\cdot)})^{(s)}}=\sup_{\substack{f\in L^{sp(\cdot)}:\\ f\ne0}}
\frac{\|Mf\|_{(L^{p(\cdot)})^{(s)}}}{\|f\|_{(L^{p(\cdot)})^{(s)}}}.
\]
However, the equivalence of the ``convexified'' quasi-norm $\|\cdot\|_{(L^{p(\cdot)})^{(s)}}$ and the
Luxemburg-Nakano quasi-norm $\|\cdot\|_{sp(\cdot)}$, as established in Theorem~\ref{th:equiv-QN}, implies that
\[
c(p_-,s)^{-1}\|M\|_{L^{sp(\cdot)}\to L^{sp(\cdot)}} \le \|M\|_{(L^{p(\cdot)})^{(s)}\to(L^{p(\cdot)})^{(s)}} \le
c(p_-,s)\|M\|_{L^{sp(\cdot)}\to L^{sp(\cdot)}},
\]
where 
\[
c(p_-,s)=2^{
\max\{1/sp_-,1\}+(1/s)\max\{1/p_-,1\}
}.
\]
Due to this relation, the boundedness of $M$ in the sense of the ``convexified'' quasi-norm is equivalent to that defined 
through the usual quasi-norm of the variable Lebesgue space, and hence the claim of the corollary follows.
\end{proof}

As we noted at the end of Subsection~\ref{subsec:example}, one may alternatively consider the maximum quasi-norm 
$\|\cdot\|_{p(\cdot)}^\text{max}$ on the spaces $L^{p(\cdot)}(\Omega,d,\mu)$ satisfying $p_->0$, with respect to 
which these spaces are again quasi-Banach lattices in possession of the Fatou property. Thus, we apply the main 
Theorem~\ref{main-theorem} to a variable Lebesgue space with the maximum quasi-norm and immediately get a version 
of Corollary~\ref{cor:L^p} for the $\|\cdot\|_{p(\cdot)}^\text{max}$ implied---in this case, the translation from the language 
of convexifications is automatic due to the equality of $\|\cdot\|_{p(\cdot)}^\text{max}$ and the ``convexified'' quasi-norm, 
established in Lemma~\ref{le:equal-QN}. 

The above two applications of the main result to variable Lebesgue spaces with different underlying quasi-norms complete
this work.

\subsection*{Acknowledgements}
This work is funded by national funds through the FCT -- Funda\c{c}\~ao para a Ci\^encia e a Tecnologia, I.P., 
under the scope of the projects UIDB/00297/2020 (\url{https://doi.org/10.54499/UIDB/00297/2020}) and 
UIDP/00297/2020 (\url{https://doi.org/10.54499/UIDP/00297/2020}) (Center for Mathematics and Applications)
and  under the scope of the PhD scholarship UI/BD/154284/2022.

I would like to thank Dr.~Oleksiy Karlovych, who supervised this work as a part of my PhD research, for 
encouraging me to address the self-improvement problem and our insightful conversations all the way through.

\bibliographystyle{elsarticle-harv} 
\bibliography{AS-bibliography}

\begin{thebibliography}{27}
\expandafter\ifx\csname natexlab\endcsname\relax\def\natexlab#1{#1}\fi
\providecommand{\url}[1]{\texttt{#1}}
\providecommand{\href}[2]{#2}
\providecommand{\path}[1]{#1}
\providecommand{\DOIprefix}{doi:}
\providecommand{\ArXivprefix}{arXiv:}
\providecommand{\URLprefix}{URL: }
\providecommand{\Pubmedprefix}{pmid:}
\providecommand{\doi}[1]{\href{http://dx.doi.org/#1}{\path{#1}}}
\providecommand{\Pubmed}[1]{\href{pmid:#1}{\path{#1}}}
\providecommand{\bibinfo}[2]{#2}
\ifx\xfnm\relax \def\xfnm[#1]{\unskip,\space#1}\fi
\bibitem[{Anderson et~al.(2017)Anderson, Hyt\"{o}nen and Tapiola}]{AHT17}
\bibinfo{author}{Anderson, T.C.}, \bibinfo{author}{Hyt\"{o}nen, T.},
  \bibinfo{author}{Tapiola, O.}, \bibinfo{year}{2017}.
\newblock \bibinfo{title}{Weak {$A_\infty$} weights and weak reverse
  {H}\"{o}lder property in a space of homogeneous type}.
\newblock \bibinfo{journal}{J. Geom. Anal.} \bibinfo{volume}{27},
  \bibinfo{pages}{95--119}.
\newblock \DOIprefix\doi{10.1007/s12220-015-9675-6}.
\bibitem[{Castillo and Rafeiro(2016)}]{CR16}
\bibinfo{author}{Castillo, R.E.}, \bibinfo{author}{Rafeiro, H.},
  \bibinfo{year}{2016}.
\newblock \bibinfo{title}{An introductory course in {L}ebesgue spaces}.
\newblock CMS Books in Mathematics/Ouvrages de Math\'{e}matiques de la SMC,
  \bibinfo{publisher}{Springer, [Cham]}.
\newblock \DOIprefix\doi{10.1007/978-3-319-30034-4}.
\bibitem[{Christ(1990a)}]{C90b}
\bibinfo{author}{Christ, M.}, \bibinfo{year}{1990}a.
\newblock \bibinfo{title}{Lectures on singular integral operators}.
  volume~\bibinfo{volume}{77} of \textit{\bibinfo{series}{CBMS Regional
  Conference Series in Mathematics}}.
\newblock \bibinfo{publisher}{Published for the Conference Board of the
  Mathematical Sciences, Washington, DC; by the American Mathematical Society,
  Providence, RI}.
\bibitem[{Christ(1990b)}]{C90a}
\bibinfo{author}{Christ, M.}, \bibinfo{year}{1990}b.
\newblock \bibinfo{title}{A {$T(b)$} theorem with remarks on analytic capacity
  and the {C}auchy integral}.
\newblock \bibinfo{journal}{Colloq. Math.} \bibinfo{volume}{60/61},
  \bibinfo{pages}{601--628}.
\newblock \DOIprefix\doi{10.4064/cm-60-61-2-601-628}.
\bibitem[{Coifman and Weiss(1971)}]{CW71}
\bibinfo{author}{Coifman, R.R.}, \bibinfo{author}{Weiss, G.},
  \bibinfo{year}{1971}.
\newblock \bibinfo{title}{Analyse harmonique non-commutative sur certains
  espaces homog\`enes}. volume \bibinfo{volume}{Vol. 242} of
  \textit{\bibinfo{series}{Lecture Notes in Mathematics}}.
\newblock \bibinfo{publisher}{Springer-Verlag, Berlin-New York}.
\newblock \bibinfo{note}{\'{E}tude de certaines int\'{e}grales singuli\`eres}.
\bibitem[{Coifman and Weiss(1977)}]{CW77}
\bibinfo{author}{Coifman, R.R.}, \bibinfo{author}{Weiss, G.},
  \bibinfo{year}{1977}.
\newblock \bibinfo{title}{Extensions of {H}ardy spaces and their use in
  analysis}.
\newblock \bibinfo{journal}{Bull. Amer. Math. Soc.} \bibinfo{volume}{83},
  \bibinfo{pages}{569--645}.
\newblock \DOIprefix\doi{10.1090/S0002-9904-1977-14325-5}.
\bibitem[{Cruz-Uribe and Fiorenza(2013)}]{CF13}
\bibinfo{author}{Cruz-Uribe, D.V.}, \bibinfo{author}{Fiorenza, A.},
  \bibinfo{year}{2013}.
\newblock \bibinfo{title}{Variable {L}ebesgue spaces}.
\newblock Applied and Numerical Harmonic Analysis,
  \bibinfo{publisher}{Birkh\"{a}user/Springer, Heidelberg}.
\newblock \DOIprefix\doi{10.1007/978-3-0348-0548-3}.
  \bibinfo{note}{{F}oundations and harmonic analysis}.
\bibitem[{Cruz-Uribe et~al.(2011)Cruz-Uribe, Martell and P\'{e}rez}]{CMP11}
\bibinfo{author}{Cruz-Uribe, D.V.}, \bibinfo{author}{Martell, J.M.},
  \bibinfo{author}{P\'{e}rez, C.}, \bibinfo{year}{2011}.
\newblock \bibinfo{title}{Weights, extrapolation and the theory of {R}ubio de
  {F}rancia}. volume \bibinfo{volume}{215} of \textit{\bibinfo{series}{Operator
  Theory: Advances and Applications}}.
\newblock \bibinfo{publisher}{Birkh\"{a}user/Springer Basel AG, Basel}.
\newblock \DOIprefix\doi{10.1007/978-3-0348-0072-3}.
\bibitem[{Cruz-Uribe and Shukla(2018)}]{CS18}
\bibinfo{author}{Cruz-Uribe, D.V.}, \bibinfo{author}{Shukla, P.},
  \bibinfo{year}{2018}.
\newblock \bibinfo{title}{The boundedness of fractional maximal operators on
  variable {L}ebesgue spaces over spaces of homogeneous type}.
\newblock \bibinfo{journal}{Studia Math.} \bibinfo{volume}{242},
  \bibinfo{pages}{109--139}.
\newblock \DOIprefix\doi{10.4064/sm8556-6-2017}.
\bibitem[{Diening(2005)}]{D05}
\bibinfo{author}{Diening, L.}, \bibinfo{year}{2005}.
\newblock \bibinfo{title}{Maximal function on {M}usielak-{O}rlicz spaces and
  generalized {L}ebesgue spaces}.
\newblock \bibinfo{journal}{Bull. Sci. Math.} \bibinfo{volume}{129},
  \bibinfo{pages}{657--700}.
\newblock \DOIprefix\doi{10.1016/j.bulsci.2003.10.003}.
\bibitem[{Diening et~al.(2011)Diening, Harjulehto, H\"{a}st\"{o} and
  R\r{u}\v{z}i\v{c}ka}]{DHHR11}
\bibinfo{author}{Diening, L.}, \bibinfo{author}{Harjulehto, P.},
  \bibinfo{author}{H\"{a}st\"{o}, P.}, \bibinfo{author}{R\r{u}\v{z}i\v{c}ka,
  M.}, \bibinfo{year}{2011}.
\newblock \bibinfo{title}{Lebesgue and {S}obolev spaces with variable
  exponents}. volume \bibinfo{volume}{2017} of \textit{\bibinfo{series}{Lecture
  Notes in Mathematics}}.
\newblock \bibinfo{publisher}{Springer, Heidelberg}.
\newblock \DOIprefix\doi{10.1007/978-3-642-18363-8}.
\bibitem[{Folland(1999)}]{F99}
\bibinfo{author}{Folland, G.B.}, \bibinfo{year}{1999}.
\newblock \bibinfo{title}{Real Analysis: Modern Techniques and Their
  Applications}.
\newblock \bibinfo{publisher}{Wiley}, \bibinfo{address}{New York}.
\bibitem[{Hyt\"{o}nen and Kairema(2012)}]{HK12}
\bibinfo{author}{Hyt\"{o}nen, T.}, \bibinfo{author}{Kairema, A.},
  \bibinfo{year}{2012}.
\newblock \bibinfo{title}{Systems of dyadic cubes in a doubling metric space}.
\newblock \bibinfo{journal}{Colloq. Math.} \bibinfo{volume}{126},
  \bibinfo{pages}{1--33}.
\newblock \DOIprefix\doi{10.4064/cm126-1-1}.
\bibitem[{Kalton et~al.(1984)Kalton, Peck and Roberts}]{KPR84}
\bibinfo{author}{Kalton, N.J.}, \bibinfo{author}{Peck, N.T.},
  \bibinfo{author}{Roberts, J.W.}, \bibinfo{year}{1984}.
\newblock \bibinfo{title}{An {$F$}-space sampler}. volume~\bibinfo{volume}{89}
  of \textit{\bibinfo{series}{London Mathematical Society Lecture Note
  Series}}.
\newblock \bibinfo{publisher}{Cambridge University Press, Cambridge}.
\newblock \DOIprefix\doi{10.1017/CBO9780511662447}.
\bibitem[{Kantorovich and Akilov(1982)}]{KA82}
\bibinfo{author}{Kantorovich, L.V.}, \bibinfo{author}{Akilov, G.P.},
  \bibinfo{year}{1982}.
\newblock \bibinfo{title}{Functional analysis}.
\newblock \bibinfo{edition}{Second} ed., \bibinfo{publisher}{Pergamon Press,
  Oxford-Elmsford, N.Y.}
\newblock \bibinfo{note}{Translated from the Russian by Howard L. Silcock}.
\bibitem[{Karlovich(2019)}]{K19}
\bibinfo{author}{Karlovich, A.}, \bibinfo{year}{2019}.
\newblock \bibinfo{title}{Hardy-{L}ittlewood maximal operator on the associate
  space of a {B}anach function space}.
\newblock \bibinfo{journal}{Real Anal. Exchange} \bibinfo{volume}{44},
  \bibinfo{pages}{119--140}.
\newblock \DOIprefix\doi{10.14321/realanalexch.44.1.0119}.
\bibitem[{Kempka and Vyb\'{\i}ral(2014)}]{KV14}
\bibinfo{author}{Kempka, H.}, \bibinfo{author}{Vyb\'{\i}ral, J.},
  \bibinfo{year}{2014}.
\newblock \bibinfo{title}{Lorentz spaces with variable exponents}.
\newblock \bibinfo{journal}{Math. Nachr.} \bibinfo{volume}{287},
  \bibinfo{pages}{938--954}.
\newblock \DOIprefix\doi{10.1002/mana.201200278}.
\bibitem[{König(1986)}]{K86}
\bibinfo{author}{König, H.}, \bibinfo{year}{1986}.
\newblock \bibinfo{title}{Eigenvalue distribution of compact operators}.
  volume~\bibinfo{volume}{16} of \textit{\bibinfo{series}{Operator Theory:
  Advances and Applications}}.
\newblock \bibinfo{publisher}{Birkh\"{a}user Verlag, Basel}.
\newblock \DOIprefix\doi{10.1007/978-3-0348-6278-3}.
\bibitem[{Lerner(2024)}]{L24}
\bibinfo{author}{Lerner, A.K.}, \bibinfo{year}{2024}.
\newblock \bibinfo{title}{A boundedness criterion for the maximal operator on
  variable {L}ebesgue spaces}.
\newblock \bibinfo{journal}{J. Anal. Math.} \bibinfo{volume}{accepted}.
\bibitem[{Lerner and Ombrosi(2010)}]{LO10}
\bibinfo{author}{Lerner, A.K.}, \bibinfo{author}{Ombrosi, S.},
  \bibinfo{year}{2010}.
\newblock \bibinfo{title}{A boundedness criterion for general maximal
  operators}.
\newblock \bibinfo{journal}{Publ. Mat.} \bibinfo{volume}{54},
  \bibinfo{pages}{53--71}.
\newblock \DOIprefix\doi{10.5565/PUBLMAT\_54110\_03}.
\bibitem[{Lerner and P\'{e}rez(2007)}]{LP07}
\bibinfo{author}{Lerner, A.K.}, \bibinfo{author}{P\'{e}rez, C.},
  \bibinfo{year}{2007}.
\newblock \bibinfo{title}{A new characterization of the {M}uckenhoupt {$A_p$}
  weights through an extension of the {L}orentz-{S}himogaki theorem}.
\newblock \bibinfo{journal}{Indiana Univ. Math. J.} \bibinfo{volume}{56},
  \bibinfo{pages}{2697--2722}.
\newblock \DOIprefix\doi{10.1512/iumj.2007.56.3112}.
\bibitem[{Lorist and Nieraeth(2024a)}]{LN23}
\bibinfo{author}{Lorist, E.}, \bibinfo{author}{Nieraeth, Z.},
  \bibinfo{year}{2024}a.
\newblock \bibinfo{title}{Banach function spaces done right}.
\newblock \bibinfo{journal}{Indag. Math. (N.S.)} \bibinfo{volume}{35},
  \bibinfo{pages}{247--268}.
\newblock \DOIprefix\doi{10.1016/j.indag.2023.11.004}.
\bibitem[{Lorist and Nieraeth(2024b)}]{LN24}
\bibinfo{author}{Lorist, E.}, \bibinfo{author}{Nieraeth, Z.},
  \bibinfo{year}{2024}b.
\newblock \bibinfo{title}{Extrapolation of compactness on banach function
  spaces}.
\newblock \bibinfo{journal}{Journal of Fourier Analysis and Applications}
  \bibinfo{volume}{30}, \bibinfo{pages}{Paper No. 30}.
\newblock \DOIprefix\doi{10.1007/s00041-024-10087-x}.
\bibitem[{Mac\'{\i}as and Segovia(1979)}]{MS79}
\bibinfo{author}{Mac\'{\i}as, R.A.}, \bibinfo{author}{Segovia, C.},
  \bibinfo{year}{1979}.
\newblock \bibinfo{title}{Lipschitz functions on spaces of homogeneous type}.
\newblock \bibinfo{journal}{Adv. in Math.} \bibinfo{volume}{33},
  \bibinfo{pages}{257--270}.
\newblock \DOIprefix\doi{10.1016/0001-8708(79)90012-4}.
\bibitem[{Maligranda(2004)}]{M04}
\bibinfo{author}{Maligranda, L.}, \bibinfo{year}{2004}.
\newblock \bibinfo{title}{Type, cotype and convexity properties of
  quasi-{B}anach spaces}, in: \bibinfo{booktitle}{Banach and function spaces}.
  \bibinfo{publisher}{Yokohama Publ., Yokohama}, pp. \bibinfo{pages}{83--120}.
\bibitem[{Nieraeth(2023)}]{N23}
\bibinfo{author}{Nieraeth, Z.}, \bibinfo{year}{2023}.
\newblock \bibinfo{title}{Extrapolation in general quasi-{B}anach function
  spaces}.
\newblock \bibinfo{journal}{J. Funct. Anal.} \bibinfo{volume}{285},
  \bibinfo{pages}{Paper No. 110130, 109}.
\newblock \DOIprefix\doi{10.1016/j.jfa.2023.110130}.
\bibitem[{Paluszy\'{n}ski and Stempak(2009)}]{PS09}
\bibinfo{author}{Paluszy\'{n}ski, M.}, \bibinfo{author}{Stempak, K.},
  \bibinfo{year}{2009}.
\newblock \bibinfo{title}{On quasi-metric and metric spaces}.
\newblock \bibinfo{journal}{Proc. Amer. Math. Soc.} \bibinfo{volume}{137},
  \bibinfo{pages}{4307--4312}.
\newblock \DOIprefix\doi{10.1090/S0002-9939-09-10058-8}.

\end{thebibliography}
\end{document}